\documentclass[reqno, oneside, 11pt]{amsart}

\usepackage[margin=1.18in]{geometry}

\usepackage{color}
\usepackage{amsthm,amsmath,amssymb,mathrsfs,eucal}
\usepackage[dvips,colorlinks,hypertex]{hyperref}

\newtheorem{thm}{Theorem}[section]

\newtheorem{lem}[thm]{Lemma}
\newtheorem{prop}[thm]{Proposition}
\theoremstyle{definition}

\newtheorem{aspt}[thm]{Assumption}
\theoremstyle{remark}
\newtheorem{rem}[thm]{Remark}

\begin{document}

\title{Discrete time nonlinear filters with informative observations
are stable}
\author{Ramon van Handel}
\address{ORFE Department\\
Princeton University\\
Princeton, NJ 08544\\
USA}
\email{rvan@princeton.edu}

\subjclass[2000]{Primary        93E11;  
                secondary       60J05,
                                62M20,  
                                93E15}  

\keywords{nonlinear filtering, prediction, asymptotic stability,
hidden Markov models}

\begin{abstract}
The nonlinear filter associated with the discrete time signal-observation
model $(X_k,Y_k)$ is known to forget its initial condition as $k\to\infty$
regardless of the observation structure when the signal possesses
sufficiently strong ergodic properties.  Conversely, it stands to reason
that if the observations are sufficiently informative, then the nonlinear
filter should forget its initial condition regardless of any properties of
the signal.  We show that for observations of additive type
$Y_k=h(X_k)+\xi_k$ with invertible observation function $h$ (under mild
regularity assumptions on $h$ and on the distribution of the noise
$\xi_k$), the filter is indeed stable in a weak sense without any
assumptions at all on the signal process.  If the signal satisfies a
uniform continuity assumption, weak stability can be strengthened to
stability in total variation.
\end{abstract}

\maketitle

\section{Introduction}

Let $(E,\mathcal{B}(E))$ and $(F,\mathcal{B}(F))$ be Polish spaces endowed
with their Borel $\sigma$-fields, and let
$P:E\times\mathcal{B}(E)\to[0,1]$ be a given transition probability
kernel.  On the sequence space $\Omega = E^{\mathbb{Z}_+}\times
F^{\mathbb{Z}_+}$ with the canonical coordinate projections
$X_n(x,y)=x(n)$, $Y_n(x,y)=y(n)$, we define the family of probability
measures $\mathbf{P}^\mu$ (for any probability measure $\mu$ on $E$) such
that $(X_n)_{n\ge 0}$ is a Markov chain with initial measure $X_0\sim\mu$
and transition probability $P$, and such that $Y_n=H(X_n,\xi_n)$ for every
$n\ge 0$ where $\xi_n$ is an i.i.d.\ sequence independent of $(X_n)_{n\ge
0}$.  A time series model of this type, called a \emph{hidden Markov
model}, has a wide variety of applications in science, engineering,
statistics and finance; see, e.g., \cite{CMR05}.  The process $(X_n)_{n\ge
0}$ is known as the \emph{signal process} (and $E$ is the signal state
space), while $(Y_n)_{n\ge 0}$ is called the \emph{observation process}
(and $F$ is the observation state space).

As $E$ is Polish, we may define for every $\mu$ the regular conditional
probabilities
$$
        \pi^\mu_{n-}(\,\cdot\,):= \mathbf{P}^\mu(X_{n}\in\,\cdot\,|
        Y_0,\ldots,Y_{n-1}),\qquad
        n\ge 1,
$$
and
$$
        \pi^\mu_n(\,\cdot\,):= \mathbf{P}^\mu(X_{n}\in\,\cdot\,|
        Y_0,\ldots,Y_{n}),\qquad
        n\ge 0.
$$
Here $\pi^\mu_{n-}$ is called the \textit{one step predictor} of the
signal given the observations, while $\pi_n^\mu$ is known as the 
\textit{nonlinear filter}.  These objects play a central role in the 
statistical theory of hidden Markov models.  A question which has 
generated considerable interest in recent years is whether, as 
$n\to\infty$, the filter $\pi_n^\mu$ becomes insensitive to the choice of 
the initial measure $\mu$.  Broadly speaking, the filter is said to be 
\emph{stable} if $\pi_n^\mu$ and $\pi_n^{\nu}$ converge towards one 
another in a suitably chosen manner as $n\to\infty$ (e.g.,
$\|\pi_n^\mu-\pi_n^\nu\|_{\rm TV}\to 0$ $\mathbf{P}^\mu$-a.s.) for a large 
class of initial measures $\mu,\nu$ (e.g., for all $\mu\ll\nu$).  

The filter stability property is of significant practical interest, as the
initial measure (a Bayesian prior) may be difficult to characterize.  
When the filter is stable, we can guarantee that it will nonetheless
generate optimal estimates of the signal process after an initial
transient.  Moreover, the stability property also plays a key role in
various important auxiliary problems, such as proving consistency of
maximum likelihood estimates and proving uniform convergence of
approximate filtering algorithms.  On the other hand, the filter stability
problem poses a set of interesting mathematical questions in the theory of
nonlinear estimation, many of which have yet to be fully resolved.  An 
overview of the state-of-the-art can be found in \cite{CR09}.

Intuitively one expects filter stability to be caused by two
separate mechanisms:
\begin{enumerate}
\item If the signal process itself becomes insensitive to its initial 
condition after a long time interval (i.e., the signal is ergodic) one 
would expect the filter to inherit this property.
\item If the observations are informative, one would expect that 
the information in the observations will eventually obsolete the prior 
information contained in the initial measure.
\end{enumerate}
In the two special cases where a detailed characterization of filter 
stability is available---for linear Gaussian models \cite{OP96} and for 
finite state signals \cite{Van08a}---the notion of \emph{detectability} 
embodies precisely this intuition.  It thus seems reasonable to 
conjecture that it is true in great generality that these two mechanisms 
conspire to bring about the stability of the filter.  To date, however,
detectability conditions for filter stability are only known in the 
abovementioned special cases.  To gain further insight, it is therefore 
instructive to study each of the mechanisms separately in a general 
setting.  In particular, the two extreme cases lead to the following 
fundamental problems: (i) can we find conditions on the signal process 
such that the filter is stable regardless of the observation structure?
and (ii) can we find conditions on the observation structure such that the 
filter is stable regardless of any properties of the signal?

Various solutions to Problem (i) can be found in the literature.  It was
shown by Atar and Zeitouni \cite{AZ97} and by Del Moral and Guionnet
\cite{DmG01} that the filter is stable whenever the signal possesses a
certain strong \emph{mixing condition}, regardless of the observation
structure.  The mixing condition was weakened to some extent by Chigansky 
and Liptser \cite{ChL04}.  Somewhat surprisingly, assuming only ergodicity 
of the signal is not sufficient to guarantee stability (see \cite[section
5]{ChL04}); both the mixing condition and the condition of Chigansky and
Liptser are strictly stronger than ergodicity.  Under a mild nondegeneracy
assumption on the observations, however, ergodicity of the signal is
already sufficient to ensure stability of the filter \cite{Van08c}.

In contrast to the first problem, a solution to Problem (ii) has hitherto
been elusive.  Unlike results based on ergodicity or mixing, stability
results based on the structure of the observations have appeared only
recently in \cite{ChL06,Van08a,Van08b}. It appears, however, that it is
more natural in this context to study stability of the predictor than
stability of the filter.  In particular, the following general result was
established in \cite[proposition 3.11]{Van08b} for \emph{additive} 
observations of the form $Y_n=h(X_n) + \xi_n$, where $E=F=\mathbb{R}^n$,
$h:\mathbb{R}^n\to\mathbb{R}^n$ is a given observation function, and
$\xi_n$ is a sequence of i.i.d.\ $\mathbb{R}^n$-valued random variables
independent of $(X_n)_{n\ge 0}$.

\begin{prop}[\cite{Van08b}]
\label{prop:pred}
Suppose that the following hold:
\begin{enumerate}
\item $h$ possesses a uniformly continuous inverse; and
\item the characteristic function of $\xi_0$ vanishes nowhere.
\end{enumerate}
Then $\|\pi_{n-}^\mu-\pi_{n-}^\nu\|_{\rm BL}\xrightarrow{n\to\infty}0$
$\mathbf{P}^\mu$-a.s.\ whenever
$\mathbf{P}^\mu|_{\sigma\{(Y_k)_{k\ge 0}\}}\ll
\mathbf{P}^\nu|_{\sigma\{(Y_k)_{k\ge 0}\}}$.
\end{prop}

Here $\|\,\cdot\,\|_{\rm BL}$ denotes the dual bounded-Lipschitz distance
(to be defined below).  The assumptions of this result certainly conform
to the idea that the observations are `informative': $Y_k$ is simply a
noisy and distorted version of $X_k$.  Note also that this result places
no conditions whatsoever on the signal process $X_n$ except the Markov
property. However, the result is a statement about the one step predictor
and not about the filter.

The main purpose of this note is to point out that under the mild
additional assumption that the law of the noise variables $\xi_k$ has a
density, a slightly weaker version of proposition \ref{prop:pred} holds
also when the predictor is replaced by the filter.  We therefore provide
an affirmative answer to the conjecture that there exists a solution to
Problem (ii) above.  The proof of this result adapts a coupling argument
due to Ocone and Pardoux \cite[lemma 3.6]{OP96}.

\begin{rem}
In the continuous time setting, it is known that a result along the lines 
of proposition \ref{prop:pred} holds for the filter when the signal state 
space $E$ is assumed to be compact (under the mild assumption that the 
signal is Feller), see \cite{Van08a}.  This is not a satisfactory solution 
to Problem (ii), however, as unstable signals are ruled out in a compact 
state space.
\end{rem}

Even if one is willing to make assumptions on the signal process, the case
of non-ergodic signals has received comparatively little attention in the
literature.  Previous results in the non-ergodic setting show that the
filter is stable in the total variation distance, but only under strong
assumptions on both the signal and the observation process
\cite{BO99,OR05,HC08}.  In particular, these results only hold when the
signal to noise ratio of the observations is sufficiently large (this
appears to be an inherent restriction of the method of proof used in these
papers).  In addition to our main result, we will show that the filter is
stable in the total variation distance under significantly weaker
assumptions than have been required in previous work.  In particular, our
result holds for an arbitrary signal to noise ratio.  To this end, we must
investigate when the convergence in the dual bounded-Lipschitz distance in
our main result can be strengthened to convergence in the total variation
distance. For this purpose we will introduce a suitable uniform continuity
assumption on the signal transition kernel.

\section{Notation and Main Results}

\subsection{Notation}
\label{sec:notation}

For any Polish space $S$ endowed with a complete metric $d$
(when $S=\mathbb{R}^n$ we will always choose the Euclidean metric
$d(x,y)=\|x-y\|$), define
$$
        \|f\|_\infty = \sup_x|f(x)|,\qquad
        \|f\|_L = \sup_{x\ne y}\frac{|f(x)-f(y)|}{d(x,y)}
	\qquad
	\mbox{for any }f:S\to\mathbb{R}.
$$
If $\|f\|_L<\infty$, the function $f$ is Lipschitz continuous.  Denote by
$\mathrm{Lip}=\{f:S\to\mathbb{R}:\|f\|_\infty\le 1\mbox{ and }\|f\|_L\le 
1\}$ the unit ball in the space of $1$-Lipschitz functions.  Then for any 
two probability measures $\mu,\nu$ on $S$, the dual bounded-Lipschitz norm 
is defined as
$$
	\|\mu-\nu\|_{\rm BL} = \sup_{f\in\mathrm{Lip}}\left|
		\int f(x)\,\mu(dx) - \int f(x)\,\nu(dx)
	\right|. 
$$ 
The supremum can equivalently be taken over a countable subfamily 
$\mathrm{Lip}_0\subset\mathrm{Lip}$ ($\mathrm{Lip}_0$ does not depend on 
$\mu,\nu$) \cite[lemma A.1]{Van08b}.  As usual, the total variation norm 
is defined as
$$
	\|\mu-\nu\|_{\rm TV} = \sup_{\|f\|_\infty\le 1}\left|
		\int f(x)\,\mu(dx) - \int f(x)\,\nu(dx)
	\right|.		
$$
Also in this case the supremum can be replaced by the supremum over a 
countable subfamily $B_0\subset\{f:\|f\|_\infty\le 1\}$ (along the lines 
of \cite[lemma 4.1]{Van08b}).

\subsection{Main Results}

In the following, we will work with a hidden Markov model where the 
signal state space $E$ is a Polish space with complete metric $d$, 
and the observation state space is Euclidean $F=\mathbb{R}^n$.  We 
consider additive observations of the form $Y_k=h(X_k)+\xi_k$ for all 
$k\ge 0$, where $h:E\to\mathbb{R}^n$ is the observation function and
$\xi_k$ is a sequence of i.i.d.\ $\mathbb{R}^n$-valued random variables
which are independent of the signal.  The signal transition kernel $P$ is 
fixed at the outset and is not presumed to satisfy any assumptions until 
further notice.

In our main result, we will impose the following assumption.

\begin{aspt}
\label{aspt}
The following hold:
\begin{enumerate}
\item The observation function $h$ possesses a uniformly continuous 
inverse.
\item The law of $\xi_k$ has a density $q_\xi$ with respect to the
Lebesgue measure on $\mathbb{R}^n$.
\item The Fourier transform of $q_\xi$ vanishes nowhere.
\end{enumerate}
\end{aspt}

Note that this is an assumption on the observations only: nothing at all
is assumed about the signal at this point. Our main result thus
holds regardless of any properties of the signal.

\begin{thm}
\label{thm:main}
Suppose that assumption \ref{aspt} holds.  Then
$$
	\mathbf{E}^\mu(\|\pi_{n}^\mu-\pi_{n}^\nu\|_{\rm BL})
	\xrightarrow{n\to\infty}0
	\qquad\mbox{whenever}
	\qquad
	\mathbf{P}^\mu|_{\sigma\{(Y_k)_{k\ge 0}\}}\ll
	\mathbf{P}^\nu|_{\sigma\{(Y_k)_{k\ge 0}\}}.
$$
\end{thm}

In order to strengthen the convergence to the total variation distance, we
do need to impose an assumption on the signal. The following is
essentially a uniform strong Feller assumption.

\begin{aspt}
\label{aspt2}
The signal transition kernel $P$ satisfies
$$
	\|P(x_n,\cdot)-P(y_n,\cdot)\|_{\rm TV}\xrightarrow{n\to\infty}0
	\qquad\mbox{whenever}\qquad
	d(x_n,y_n)\xrightarrow{n\to\infty}0.
$$
In other words, the measure-valued map $x\mapsto P(x,\cdot)$ is uniformly
continuous for the total variation distance on the space of probability 
measures.
\end{aspt}

We now obtain the following result.

\begin{thm}
\label{thm:tv}
Suppose that assumptions \ref{aspt} and \ref{aspt2} hold.  Then
$$
	\|\pi_{n}^\mu-\pi_{n}^\nu\|_{\rm TV}
	\xrightarrow{n\to\infty}0\quad\mathbf{P}^\mu\mbox{-a.s.}
	\qquad\mbox{whenever}
	\qquad
	\mathbf{P}^\mu|_{\sigma\{(Y_k)_{k\ge 0}\}}\ll
	\mathbf{P}^\nu|_{\sigma\{(Y_k)_{k\ge 0}\}}.
$$
\end{thm}

A typical example where assumption \ref{aspt2} holds is the following.

\begin{prop}
\label{prop:ar}
Let $E=\mathbb{R}^m$. Suppose that the signal is defined by the 
recursion
$$
	X_{k+1} = b(X_k) + \sigma(X_k)\,\eta_k,
$$
where $b,\sigma,\eta_k$ satisfy the following assumptions:
\begin{enumerate}
\item $b:\mathbb{R}^m\to\mathbb{R}^m$ and 
$\sigma:\mathbb{R}^m\to\mathbb{R}^{m\times m}$ are uniformly continuous.
\item $\sigma$ is uniformly bounded from below:
$\|\sigma(x)v\|\ge \alpha\|v\|$ for all $x,v\in\mathbb{R}^m$ and some
$\alpha>0$.
\item $\eta_k$ are i.i.d.\ $\mathbb{R}^m$-valued random variables, whose 
law possesses a density $q_\eta$ with respect to the Lebesgue measure on 
$\mathbb{R}^m$.
\end{enumerate}
Then assumption \ref{aspt2} holds.  If $q_\eta$ and $q_\xi$ are
strictly positive and assumption \ref{aspt} holds, then
$$
	\|\pi_{n}^\mu-\pi_{n}^\nu\|_{\rm TV}
	\xrightarrow{n\to\infty}0\quad\mathbf{P}^\gamma\mbox{-a.s.}
	\qquad\mbox{for any }
	\mu,\nu,\gamma.
$$
\end{prop}

This result should be compared to the main results in
\cite{BO99,OR05,HC08}, where a very similar model is investigated.
However, in these references total variation stability is proved only when
the signal to noise ratio is sufficiently high.  This is an artefact of
the quantitative method of proof where two rates of expansion are 
compared: the filter is stable if one of the rates `wins', which leads to 
a requirement on the signal to noise ratio.  Our qualitative approach does 
not depend on the signal to noise ratio, however, so that evidently the 
assumptions required for the balancing of rates are stronger than is 
needed for filter stability (see \cite{Van08b} for further discussion).  
On the other hand, our approach can not provide an estimate of the rate of 
stability.

\begin{rem}
Our results require that the observation state space is Euclidean
$F=\mathbb{R}^n$, as it relies on properties of convolutions in
$\mathbb{R}^n$ (an extension to the case where $F$ is a locally compact 
abelian group may be feasible).  In contrast, we have only assumed that 
the signal state space is Polish.  Note, however, that assumption 
\ref{aspt} requires the existence of a uniformly continuous map 
$h^{-1}:\mathbb{R}^n\to E$ such that $h^{-1}(h(x))=x$ for all $x\in E$.
In particular, $h$ is an embedding of $E$ into $\mathbb{R}^n$, so that $E$ 
can not be larger (e.g., of higher dimension) than $\mathbb{R}^n$.
This is to be expected, of course, as filter stability in the case where 
$h$ is not invertible must depend on specific properties of the signal 
process such as observability or ergodicity.
\end{rem}

\begin{rem}
Theorems \ref{thm:main} and \ref{thm:tv} provide stability of the filter 
whenever the initial measures satisfy 
$\mathbf{P}^\mu|_{\sigma\{(Y_k)_{k\ge 0}\}}\ll
\mathbf{P}^\nu|_{\sigma\{(Y_k)_{k\ge 0}\}}$.  Absolute continuity of 
the initial measures $\mu\ll\nu$ is sufficient for this to hold, but is 
not always necessary.  For example, if $\mu P^k\ll\nu P^k$ for some $k>0$ 
and the observation density $q_\xi$ is strictly positive, then it is not 
difficult to prove that $\mathbf{P}^\mu|_{\sigma\{(Y_k)_{k\ge 0}\}}\ll
\mathbf{P}^\nu|_{\sigma\{(Y_k)_{k\ge 0}\}}$ also.  In particular, if the 
signal possesses a strictly positive transition density, then $\mu P\sim 
\nu P$ for every $\mu,\nu$ and we obtain stability for arbitrary initial 
measures provided $q_\xi>0$.  This is the case, for example, in the 
setting of proposition \ref{prop:ar}.
\end{rem}

\begin{rem}
Our results do not give a rate of stability, while most previous work on 
filter stability gives exponential convergence rates.  The following 
simple example demonstrates that exponential stability can not be expected 
in the general setting of this paper.

Let $E=F=\mathbb{R}$ and $Y_k=X_k+\xi_k$, where $\xi_k$ are i.i.d.\ 
$N(0,1)$ and $X_k=X_0$ for all $k$.  This setting certainly satisfies the 
requirements of theorem \ref{thm:main}.  We choose 
$\mu=N(\alpha,\sigma^2)$ and $\nu=N(\beta,\sigma^2)$ for some
$\alpha,\beta,\sigma\in\mathbb{R}$ (so $\mathbf{P}^\mu\sim\mathbf{P}^\nu$).  
Linear filtering theory shows that $\pi_k^\mu$ is a random Gaussian 
measure with mean $Z_k^\mu$ and variance $V_k^\mu$ given by
$$
	Z_k^\mu = \frac{\alpha}{1 + \sigma^2(k+1)} + 
	\frac{\sigma^2(k+1)}{1+\sigma^2(k+1)}\cdot
	\frac{1}{k+1}\sum_{\ell=0}^kY_\ell,
	\qquad\qquad
	V_k^\mu = \frac{\sigma^2}{1+\sigma^2(k+1)},
$$
and similarly for $\pi_k^\nu,Z_k^\nu,V_k^\nu$ where $\alpha$ is replaced 
by $\beta$.  Note that by the law of large numbers, the second term in the 
expression for $Z_k^\mu$ (and $Z_k^\nu$) converges $\mathbf{P}^\mu$-a.s.\ 
to $X_0$.  But
$$
	\|\pi_n^\mu-\pi_n^\nu\|_{\rm BL} \ge
	\left|\int \cos(x)\,\pi_n^\mu(dx) -
	\int \cos(x)\,\pi_n^\nu(dx)\right| =
	e^{-V_n^\mu/2}|\cos(Z_n^\mu)-\cos(Z_n^\nu)|.
$$
Noting that $\cos(Z_n^\nu) = \cos(Z_n^\mu) - 
\sin(Z_n^\mu)\,(\beta-\alpha)/(1+\sigma^2(n+1))
+ o(n^{-1})$, we find 
that
$$
	\liminf_{n\to\infty}n\|\pi_n^\mu-\pi_n^\nu\|_{\rm BL}
	\ge \liminf_{n\to\infty} n\,e^{-V_n^\mu/2}|\cos(Z_n^\mu)-\cos(Z_n^\nu)|
	= \frac{|\beta-\alpha|}{\sigma^2}\,|\sin(X_0)|>0
	\quad\mathbf{P}^\mu\mbox{-a.s.}
$$
By Fatou's lemma $\liminf_{n\to\infty}n\,\mathbf{E}^\mu(
\|\pi_n^\mu-\pi_n^\nu\|_{\rm BL})>0$, so that evidently the stability rate 
of the filter is at best of order $O(n^{-1})$ and is certainly not 
exponential.

It is interesting to note that by \cite[p.\ 528, theorem 4]{Shi96} and by 
the equivalence of the Hellinger and total variation distances,
$\sum_{n=0}^\infty \|\pi^\mu_{n-} h^{-1}\ast\xi-\pi^\nu_{n-} 
h^{-1}\ast\xi\|_{\rm TV}^2<\infty$ $\mathbf{P}^\mu$-a.s.  The convergence 
of the expression $\|\pi^\mu_{n-} h^{-1}\ast\xi-\pi^\nu_{n-} 
h^{-1}\ast\xi\|_{\rm TV}$ which appears in the proof of lemma 
\ref{lem:blackwell} below is therefore generally, in a sense, not much 
worse than $o(n^{-1/2})$.  It is unclear, however, whether this property 
survives the subsequent manipulations that lead to stability of the filter.
\end{rem}

\section{Proof of Theorem \ref{thm:main}}

Let us begin by recalling a part of the proof of proposition 
\ref{prop:pred}.  

\begin{lem}
\label{lem:blackwell}
Suppose that the characteristic function of $\xi_k$ vanishes nowhere, and 
that moreover $\mathbf{P}^\mu|_{\sigma\{(Y_k)_{k\ge 0}\}}\ll
\mathbf{P}^\nu|_{\sigma\{(Y_k)_{k\ge 0}\}}$.
Then $\|\pi_{n-}^\mu h^{-1}-\pi_{n-}^\nu h^{-1}\|_{\rm 
BL}\xrightarrow[n\to\infty]{}0$ $\mathbf{P}^\mu$-a.s.
\end{lem}

\begin{proof}
Denote the law of $\xi_k$ as $\xi$.  It is easily verified
that for any probability measure $\rho$ 
$$
        \mathbf{P}^\rho(Y_{n}\in\,\cdot\,|Y_0,\ldots,Y_{n-1})
        = \pi^\rho_{n-}h^{-1}\ast\xi,
$$
where $\ast$ denotes convolution.  A classic result of Blackwell and 
Dubins \cite[section 2]{BD62} shows that
$$
        \|\pi^\mu_{n-} h^{-1}\ast\xi-
        \pi^\nu_{n-} h^{-1}\ast\xi\|_{\rm TV}
        \xrightarrow{n\to\infty}0\quad
        \mathbf{P}^\mu\mbox{-a.s.},
$$
where $\|\,\cdot\,\|_{\rm TV}$ is the total variation norm.
The result now follows from \cite[proposition C.2]{Van08b}.
\end{proof}

To proceed, recall that due to the Bayes formula 
(e.g., \cite[section 3.2.2]{CMR05})
$$
	\int f(x)\,\pi^\mu_{n}(dx) =
	\frac{\int f(x)\,q_\xi(Y_n-h(x))\,\pi^\mu_{n-}(dx)
	}{\int q_\xi(Y_n-h(x))\,\pi^\mu_{n-}(dx)}\qquad
	\mbox{for all bounded }f:E\to\mathbb{R}
	\qquad
	\mathbf{P}^\mu\mbox{-a.s.}
$$
Note that the denominator of this expression is strictly positive 
$\mathbf{P}^\mu$-a.s.  Moreover, we have
$$
	\mathbf{E}^\mu(f(Y_n)|Y_0,\ldots,Y_{n-1}) =
	\int f(y)\,q_\xi(y-h(x))\,\pi^\mu_{n-}(dx)\,dy
	\qquad
	\mathbf{P}^\mu\mbox{-a.s.}
$$
for any bounded function $f:\mathbb{R}^n\to\mathbb{R}$ as in the proof of 
lemma \ref{lem:blackwell}.  Therefore, it follows from the disintegration 
of measures that $\mathbf{P}^\mu$-a.s.
\begin{multline*}
	\mathbf{E}^\mu(\|\pi_{n}^\mu-\pi_{n}^\nu\|_{\rm BL}|
	Y_0,\ldots,Y_{n-1}) = \\
	\int
	\sup_{f\in\mathrm{Lip}_0}\left|
	\frac{\int f(x)\,q_\xi(y-h(x))\,\pi^\mu_{n-}(dx)
	}{\int q_\xi(y-h(x))\,\pi^\mu_{n-}(dx)}-
	\frac{\int f(x)\,q_\xi(y-h(x))\,\pi^\nu_{n-}(dx)
	}{\int q_\xi(y-h(x))\,\pi^\nu_{n-}(dx)}
	\right| \\
	\times
	\left\{\int q_\xi(y-h(x))\,\pi^\mu_{n-}(dx)\right\}dy,
\end{multline*}
where it should be noted that by the assumption that
$\mathbf{P}^\mu|_{\sigma\{(Y_k)_{k\ge 0}\}}\ll
\mathbf{P}^\nu|_{\sigma\{(Y_k)_{k\ge 0}\}}$ in theorem \ref{thm:main}
(which we presume to be in force throughout) all quantities in this
expression as $\mathbf{P}^\mu$-a.s.\ uniquely defined and both
denominators are strictly positive $\mathbf{P}^\mu$-a.s.  

It will be useful for what follows to rewrite the above expression in a 
more convenient form:
\begin{multline*}
	\mathbf{E}^\mu(\|\pi_{n}^\mu-\pi_{n}^\nu\|_{\rm BL}|
	Y_0,\ldots,Y_{n-1}) = 
	\int
	\sup_{f\in\mathrm{Lip}_0}\left|
	\int f(h^{-1}(x))\,q_\xi(y-x)\,\pi^\mu_{n-}h^{-1}(dx) 
	\right. \\
	\left.\mbox{} -
	\frac{\int f(h^{-1}(x))\,q_\xi(y-x)\,\pi^\nu_{n-}h^{-1}(dx)
	}{\int q_\xi(y-x)\,\pi^\nu_{n-}h^{-1}(dx)}
	\int q_\xi(y-x)\,\pi^\mu_{n-}h^{-1}(dx)
	\right|dy.
\end{multline*}
Here we have fixed a uniformly continuous function 
$h^{-1}:\mathbb{R}^n\to E$ such that $h^{-1}(h(x))=x$ for all 
$x\in E$; the existence of this function is guaranteed by 
assumption \ref{aspt}.

\begin{lem}
\label{lem:couple}
Let $\rho,\rho'$ be two probability measures on $\mathbb{R}^n$, and let
$Z,Z'$ be (not necessarily independent) $\mathbb{R}^n$-valued random 
variables such that $Z\sim\rho$ and $Z'\sim\rho'$.  Then
\begin{multline*}
	\int\sup_{f\in\mathrm{Lip}_0}\left|
	\int f(h^{-1}(x))\,q_\xi(y-x)\,\rho(dx)-
	\frac{\int f(h^{-1}(x))\,q_\xi(y-x)\,\rho'(dx)
	}{\int q_\xi(y-x)\,\rho'(dx)}
	\int q_\xi(y-x)\,\rho(dx)
	\right| dy\\
	\le  \mathbf{E}(d(h^{-1}(Z),h^{-1}(Z'))\wedge 2) + 
	2\int\mathbf{E}(|q_\xi(y-Z)-q_\xi(y-Z')|)\,dy,
\end{multline*}
where by convention $0/0=1$.
\end{lem}

\begin{proof}
The left hand side of the expression in the statement of the lemma is
$$
	\Delta := \int\sup_{f\in\mathrm{Lip}_0}\left|
	\mathbf{E}(f(h^{-1}(Z))\,q_\xi(y-Z)) -
	\frac{\mathbf{E}(f(h^{-1}(Z'))\,q_\xi(y-Z'))}{
	\mathbf{E}(q_\xi(y-Z'))}
	\,\mathbf{E}(q_\xi(y-Z))
	\right|dy.
$$
Now note that for any $f\in\mathrm{Lip}$ and $x,y\in\mathbb{R}^n$ we
have $|f(x)-f(y)|\le d(x,y)\wedge 2$, so
\begin{multline*}
	\int\sup_{f\in\mathrm{Lip}_0}
	|\mathbf{E}(f(h^{-1}(Z))\,q_\xi(y-Z)) -
	\mathbf{E}(f(h^{-1}(Z'))\,q_\xi(y-Z))|\,dy \\ \le
	\int\mathbf{E}(\{d(h^{-1}(Z),h^{-1}(Z'))\wedge 2\}\,
		q_\xi(y-Z))\,dy
	= \mathbf{E}(d(h^{-1}(Z),h^{-1}(Z'))\wedge 2),
\end{multline*}
where we have used the Fubini-Tonelli theorem to exchange the order of 
integration.  Thus
\begin{multline*}
	\Delta \le \mathbf{E}(d(h^{-1}(Z),h^{-1}(Z'))\wedge 2) + \mbox{} \\
	\int\sup_{f\in\mathrm{Lip}_0}\left|
	\mathbf{E}(f(h^{-1}(Z'))\,q_\xi(y-Z)) -
	\frac{\mathbf{E}(f(h^{-1}(Z'))\,q_\xi(y-Z'))}{
	\mathbf{E}(q_\xi(y-Z'))}
	\,\mathbf{E}(q_\xi(y-Z))
	\right|dy.
\end{multline*}
Estimating $\mathbf{E}(f(h^{-1}(Z'))\,q_\xi(y-Z))$ by 
$\mathbf{E}(f(h^{-1}(Z'))\,q_\xi(y-Z'))$, we similarly obtain
\begin{multline*}
	\Delta \le \mathbf{E}(d(h^{-1}(Z),h^{-1}(Z'))\wedge 2) + 
	\int \mathbf{E}(|q_\xi(y-Z)-q_\xi(y-Z')|)\,dy +
	\mbox{} \\
	\int\sup_{f\in\mathrm{Lip}_0}\left|
	\mathbf{E}(f(h^{-1}(Z'))\,q_\xi(y-Z')) -
	\frac{\mathbf{E}(f(h^{-1}(Z'))\,q_\xi(y-Z'))}{
	\mathbf{E}(q_\xi(y-Z'))}
	\,\mathbf{E}(q_\xi(y-Z))
	\right|dy.
\end{multline*}
We now substitute in this expression
$$
	\mathbf{E}(f(h^{-1}(Z'))\,q_\xi(y-Z')) =
	\frac{\mathbf{E}(f(h^{-1}(Z'))\,q_\xi(y-Z'))}{
	\mathbf{E}(q_\xi(y-Z'))}
	\,\mathbf{E}(q_\xi(y-Z')),
$$
and note that $|\mathbf{E}(f(h^{-1}(Z'))\,q_\xi(y-Z'))/
\mathbf{E}(q_\xi(y-Z'))|\le 1$ whenever $\|f\|_\infty\le 1$.
The remainder of the proof is now immediate.
\end{proof}

A remarkable result due to Dudley \cite[theorem 11.7.1]{Dud02}, which 
extends the classical Skorokhod representation theorem to the 
$\|\,\cdot\,\|_{\rm BL}$-uniformity, allows us to put this lemma to good 
use.

\begin{lem}
Let $\rho_n$ and $\rho_n'$, $n\ge 0$ be two sequences of probability 
measures on $\mathbb{R}^n$ such that $\|\rho_n-\rho'_n\|_{\rm BL}\to 0$
as $n\to\infty$.  Then the following quantity
$$
	\int\!\sup_{f\in\mathrm{Lip}_0}\left|
	\int f(h^{-1}(x))\,q_\xi(y-x)\,\rho_n(dx)-
	\frac{\int f(h^{-1}(x))\,q_\xi(y-x)\,\rho_n'(dx)
	}{\int q_\xi(y-x)\,\rho_n'(dx)}
	\int q_\xi(y-x)\,\rho_n(dx)
	\right| dy
$$
converges to zero as $n\to\infty$.
\end{lem}

\begin{proof}
By \cite[theorem 11.7.1]{Dud02} we can construct two sequences of 
$\mathbb{R}^n$-valued random variables $Z_n$ and $Z_n'$, $n\ge 0$ on some 
underlying probability space such that $Z_n\sim\rho_n$ and 
$Z_n'\sim\rho_n'$ for every $n$ and $\|Z_n-Z_n'\|\to 0$ a.s.\ as 
$n\to\infty$.  By the previous lemma, the expression $\Delta_n$ in the 
statement of the present lemma is bounded by
$$
	\Delta_n\le 
	\mathbf{E}(d(h^{-1}(Z_n),h^{-1}(Z_n'))\wedge 2) + 
	2\int\mathbf{E}(|q_\xi(y-Z_n)-q_\xi(y-Z_n')|)\,dy
$$
for every $n$.  As $h^{-1}$ is uniformly continuous and $\|Z_n-Z_n'\|\to 
0$ a.s.\ as $n\to\infty$, we find that $d(h^{-1}(Z_n),h^{-1}(Z_n'))\to 0$ 
a.s.\ as $n\to\infty$.  Thus the first term evidently converges to zero by 
dominated convergence.  To deal with the second term, note that
$$
	\int |q_\xi(y-Z_n)-q_\xi(y-Z_n')|\,dy = 
	\int |q_\xi(y+Z_n'-Z_n)-q_\xi(y)|\,dy =
	\|T_{Z_n-Z_n'}q_\xi-q_\xi\|_{L^1(dy)},
$$
where $(T_zf)(x)=f(x-z)$ denotes translation.  But recall that translation 
is continuous in the $L^1$-topology \cite[proposition 8.5]{Fol99}, so we 
find that
$$
	\|T_{Z_n-Z_n'}q_\xi-q_\xi\|_{L^1(dy)}
	\xrightarrow{n\to\infty}0\quad\mbox{a.s.}
$$
On the other hand,
$$
	\|T_{Z_n-Z_n'}q_\xi-q_\xi\|_{L^1(dy)}
	\le
	2\,\|q_\xi\|_{L^1(dy)} = 2\qquad\mbox{for all }n.
$$
Dominated convergence gives
$$
	\int\mathbf{E}(|q_\xi(y-Z_n)-q_\xi(y-Z_n')|)\,dy
	= \mathbf{E}(\|T_{Z_n-Z_n'}q_\xi-q_\xi\|_{L^1(dy)})
	\xrightarrow{n\to\infty}0,
$$
where we have used the Fubini-Tonelli theorem to exchange the order
of integration.
\end{proof}

The proof of theorem \ref{thm:main} is now easily completed.  Indeed, 
under our assumptions we obtain $\|\pi_{n-}^\mu h^{-1}-\pi_{n-}^\nu 
h^{-1}\|_{\rm BL}\to 0$ $\mathbf{P}^\mu$-a.s.\ as $n\to\infty$ by 
lemma \ref{lem:blackwell}, so the previous lemma gives
$$
	\mathbf{E}^\mu(\|\pi_{n}^\mu-\pi_{n}^\nu\|_{\rm BL}|
	Y_0,\ldots,Y_{n-1}) \xrightarrow{n\to\infty} 0
	\quad
	\mathbf{P}^\mu\mbox{-a.s.}
$$
Taking the expectation with respect to $\mathbf{P}^\mu$, and noting that
$\|\pi_{n}^\mu-\pi_{n}^\nu\|_{\rm BL}\le 2$ so that the dominated 
convergence theorem applies, yields the proof.

\section{Proof of Theorem \ref{thm:tv}}

We begin by showing that the one step predictor is stable in mean total 
variation.

\begin{lem}
Suppose that assumptions \ref{aspt} and \ref{aspt2} hold.  Then
$$
        \mathbf{E}^\mu(\|\pi_{n-}^\mu-\pi_{n-}^\nu\|_{\rm TV})
        \xrightarrow{n\to\infty}0
        \qquad\mbox{whenever}
        \qquad
        \mathbf{P}^\mu|_{\sigma\{(Y_k)_{k\ge 0}\}}\ll
        \mathbf{P}^\nu|_{\sigma\{(Y_k)_{k\ge 0}\}}.
$$
\end{lem}

\begin{proof}
Recall that
$$
	\int f(x)\,\pi_{n-}^\mu(dx) =
	\int f(x')\,P(x,dx')\,\pi_{n-1}^\mu(dx)
$$
for any bounded measurable function $f$.  Therefore 
$$
	\|\pi_{n-}^\mu-\pi_{n-}^\nu\|_{\rm TV}
	=\|\pi_{n-1}^\mu-\pi_{n-1}^\mu\|_G
	:=
	\sup_{f\in G}\left|
	\int f(x)\,\pi_{n-1}^\mu(dx)-
	\int f(x)\,\pi_{n-1}^\nu(dx)
	\right|,
$$
where $G=\{Pf: f\in B_0\}$ (recall that $B_0$ is a countable family
of functions such that $\|\mu-\nu\|_{\rm TV}=\|\mu-\nu\|_{B_0}$, see
section \ref{sec:notation}).  We now claim that the family $G$ is 
uniformly bounded and uniformly equicontinuous.  Indeed, it is immediate 
that $\|f\|_\infty\le 1$ for every $f\in G$, as this is the case for 
every $f\in B_0$.  To prove uniform equicontinuity, note that
$$
	\sup_{d(x,y)\le\delta}|f(x)-f(y)|\le
	\sup_{d(x,y)\le\delta}\sup_{f\in B_0}|Pf(x)-Pf(y)|=
	\sup_{d(x,y)\le\delta}\|P(x,\cdot)-P(y,\cdot)\|_{\rm TV}:=
	\varpi_P(\delta)
$$
for every $f\in G$.  By assumption \ref{aspt2}, we evidently have 
$\varpi_P(\delta)\to 0$ as $\delta\to 0$.  Uniform equicontinuity of 
$G$ is therefore established.  To complete the proof, it remains to show 
that
$$
        \mathbf{E}^\mu(\|\pi_{n}^\mu-\pi_{n}^\nu\|_{G})
        \xrightarrow{n\to\infty}0
        \qquad\mbox{whenever}
        \qquad
        \mathbf{P}^\mu|_{\sigma\{(Y_k)_{k\ge 0}\}}\ll
        \mathbf{P}^\nu|_{\sigma\{(Y_k)_{k\ge 0}\}}.
$$
This is established precisely as in the proof of theorem 
\ref{thm:main}, however: the only modification that must be made in the 
present setting is that the term $\mathbf{E}(d(h^{-1}(Z),h^{-1}(Z'))\wedge 
2)$ in lemma \ref{lem:couple} is replaced by
$\mathbf{E}(\varpi_P(d(h^{-1}(Z),h^{-1}(Z')))\wedge 2)$.
\end{proof}

We now proceed to show stability of the filter (rather than the one step
predictor).  Stability in the mean follows trivially from the previous
lemma and the following estimate.

\begin{lem}
Whenever $\mathbf{P}^\mu|_{\sigma\{(Y_k)_{k\ge 0}\}}\ll
\mathbf{P}^\nu|_{\sigma\{(Y_k)_{k\ge 0}\}}$, we have
$$
	\mathbf{E}^\mu(\|\pi_{n}^\mu-\pi_{n}^\nu\|_{\rm TV}|Y_0,\ldots,Y_{n-1})
	\le 2\,\|\pi_{n-}^\mu-\pi_{n-}^\nu\|_{\rm TV}
	\quad\mathbf{P}^\mu\mbox{-a.s.}
$$
\end{lem}

\begin{proof}
As in the proof of theorem \ref{thm:main}, we can write
\begin{multline*}
	\mathbf{E}^\mu(\|\pi_{n}^\mu-\pi_{n}^\nu\|_{\rm TV}|
	Y_0,\ldots,Y_{n-1}) = 
	\int
	\sup_{f\in B_0}\left|
	\int f(x)\,q_\xi(y-h(x))\,\pi^\mu_{n-}(dx) 
	\right. \\
	\left.\mbox{} -
	\frac{\int f(x)\,q_\xi(y-h(x))\,\pi^\nu_{n-}(dx)
	}{\int q_\xi(y-h(x))\,\pi^\nu_{n-}(dx)}
	\int q_\xi(y-h(x))\,\pi^\mu_{n-}(dx)
	\right|dy.
\end{multline*}
It follows directly that we can estimate
$\mathbf{E}^\mu(\|\pi_{n}^\mu-\pi_{n}^\nu\|_{\rm TV}|
Y_0,\ldots,Y_{n-1})\le \Delta_1+\Delta_2$, where
\begin{equation*}
\begin{split}
	\Delta_1 &=
	\int
	\sup_{f\in B_0}\left|
	\int f(x)\,q_\xi(y-h(x))\,\pi^\mu_{n-}(dx) -
	\int f(x)\,q_\xi(y-h(x))\,\pi^\nu_{n-}(dx)
	\right|dy,\\
	\Delta_2 &=
	\int
	\sup_{f\in B_0}\left|
	\int f(x)\,q_\xi(y-h(x))\,\pi^\nu_{n-}(dx) 
	\right.\\
	&\left.\qquad\qquad\qquad\qquad
	\mbox{}-
	\frac{\int f(x)\,q_\xi(y-h(x))\,\pi^\nu_{n-}(dx)
	}{\int q_\xi(y-h(x))\,\pi^\nu_{n-}(dx)}
	\int q_\xi(y-h(x))\,\pi^\mu_{n-}(dx)
	\right|dy \\
	&\le
	\int\left|
	\int q_\xi(y-h(x))\,\pi^\mu_{n-}(dx) -
	\int q_\xi(y-h(x))\,\pi^\nu_{n-}(dx)
	\right|dy.
\end{split}
\end{equation*}
To estimate $\Delta_1$, note that
\begin{multline*}
	\sup_{f\in B_0}\left|
	\int f(x)\,q_\xi(y-h(x))\,\pi^\mu_{n-}(dx) -
	\int f(x)\,q_\xi(y-h(x))\,\pi^\nu_{n-}(dx)
	\right| \\ \le
	\sup_{f\in B_0}
	\int |f(x)|\,q_\xi(y-h(x))\,|\pi^\mu_{n-}-\pi^\nu_{n-}|(dx)
	\le
	\int q_\xi(y-h(x))\,|\pi^\mu_{n-}-\pi^\nu_{n-}|(dx).
\end{multline*}
Therefore, the Fubini-Tonelli theorem gives
$$
	\Delta_1\le
	\int \left\{\int q_\xi(y-h(x))\,dy\right\}
	|\pi^\mu_{n-}-\pi^\nu_{n-}|(dx) =
	\|\pi^\mu_{n-}-\pi^\nu_{n-}\|_{\rm TV}.
$$
$\Delta_2$ is estimated in the same fashion, and the proof is complete.
\end{proof}

We have now shown that under assumptions \ref{aspt} and \ref{aspt2},
$$
	\mathbf{E}^\mu(\|\pi_{n}^\mu-\pi_{n}^\nu\|_{\rm TV})
	\xrightarrow{n\to\infty}0
	\qquad\mbox{whenever}
	\qquad
	\mathbf{P}^\mu|_{\sigma\{(Y_k)_{k\ge 0}\}}\ll
	\mathbf{P}^\nu|_{\sigma\{(Y_k)_{k\ge 0}\}}.
$$
It remains to prove that in fact
$$
	\|\pi_{n}^\mu-\pi_{n}^\nu\|_{\rm TV}
	\xrightarrow{n\to\infty}0\quad\mathbf{P}^\mu\mbox{-a.s.}
	\qquad\mbox{whenever}
	\qquad
	\mathbf{P}^\mu|_{\sigma\{(Y_k)_{k\ge 0}\}}\ll
	\mathbf{P}^\nu|_{\sigma\{(Y_k)_{k\ge 0}\}}.
$$
Clearly it suffices to show that $\|\pi_{n}^\mu-\pi_{n}^\nu\|_{\rm TV}$
is $\mathbf{P}^\mu$-a.s.\ convergent.

To this end, let $\mu\ll\gamma$.
Then (see \cite[corollary 5.7]{Van08c})
$$
	\|\pi_{n}^\mu-\pi_{n}^\gamma\|_{\rm TV} = 
        \frac{
        \mathbf{E^\gamma}(
        |\mathbf{E^\gamma}(\frac{d\mu}{d\gamma}(X_0)\mbox{}|\mathcal{F}^Y_+\vee
                \mathcal{F}^X_{[n,\infty[})
                -      
	\mathbf{E^\gamma}(\frac{d\mu}{d\gamma}(X_0)\mbox{}|\mathcal{F}^Y_{[0,n]})|
        ~|\mathcal{F}_{[0,n]}^Y)
        }{
	\mathbf{E^\gamma}(\frac{d\mu}{d\gamma}(X_0)\mbox{}|\mathcal{F}^Y_{[0,n]})
        }
	\quad
	\mathbf{P}^\mu\mbox{-a.s.},
$$
where $\mathcal{F}^Y_{[0,n]}:=\sigma\{Y_0,\ldots,Y_n\}$,
$\mathcal{F}^Y_+:=\sigma\{Y_k:k\ge 0\}$, and $\mathcal{F}^X_{[n,\infty[}:=
\sigma\{X_k:k\ge n\}$.  When $d\mu/d\gamma$ is bounded, the numerator 
converges $\mathbf{P}^\gamma$-a.s.\ (hence $\mathbf{P}^\mu$-a.s.) by the 
martingale convergence theorem (see, e.g., \cite[theorem 2]{BD62}) while 
the denominator converges to a $\mathbf{P}^\mu$-a.s.\ strictly positive 
quantity $\mathbf{E^\gamma}(d\mu/d\gamma(X_0)\mbox{}|\mathcal{F}^Y_{[0,n]})
\to\mathbf{E^\gamma}(d\mu/d\gamma(X_0)\mbox{}|\mathcal{F}^Y_{+})>0$
$\mathbf{P}^\mu$-a.s. Evidently 
$$
	\|\pi_{n}^\mu-\pi_{n}^\gamma\|_{\rm TV}
	\xrightarrow{n\to\infty}0\quad\mathbf{P}^\mu\mbox{-a.s.}
	\qquad\mbox{whenever}\qquad\mu\ll\gamma,~\|d\mu/d\gamma\|_\infty<\infty.
$$
Now set $\gamma=(\mu+\nu)/2$, and note that
$\|d\mu/d\gamma\|_\infty\le 2$, $\|d\nu/d\gamma\|_\infty\le 2$.
Therefore
$$
	\|\pi_{n}^\mu-\pi_{n}^\gamma\|_{\rm TV}
	\xrightarrow{n\to\infty}0\quad\mathbf{P}^\mu\mbox{-a.s.},
	\qquad\quad
	\|\pi_{n}^\nu-\pi_{n}^\gamma\|_{\rm TV}
	\xrightarrow{n\to\infty}0\quad\mathbf{P}^\nu\mbox{-a.s.}
$$
If $\mathbf{P}^\mu|_{\sigma\{(Y_k)_{k\ge 0}\}}\ll
\mathbf{P}^\nu|_{\sigma\{(Y_k)_{k\ge 0}\}}$ the second statement holds 
also $\mathbf{P}^\mu$-a.s.  The proof of theorem \ref{thm:tv} is now 
easily completed by applying the triangle inequality.

\begin{rem}
From the above expression, it can be read off that under assumptions
\ref{aspt} and \ref{aspt2}
$$
	\mathbf{E}^\nu\Bigg(f(X_0)\,\Bigg|\bigcap_{n\ge 0}
	\mathcal{F}^Y_+\vee\mathcal{F}^X_{[n,\infty[}\Bigg) =
	\mathbf{E}^\nu(f(X_0)|
	\mathcal{F}^Y_+)
	\quad\mbox{whenever}\quad
	\|f\|_\infty<\infty
$$
for every $\nu$.  With a little more work, one can show that similarly
for every $\nu$
$$
	\mathbf{E}^\nu\Bigg(f(X_0,\ldots,X_k)\,\Bigg|\bigcap_{n\ge 0}
	\mathcal{F}^Y_+\vee\mathcal{F}^X_{[n,\infty[}\Bigg) =
	\mathbf{E}^\nu(f(X_0,\ldots,X_k)|
	\mathcal{F}^Y_+)
	\quad\mbox{whenever}\quad
	\|f\|_\infty<\infty,
$$
which implies that for every $\nu$
$$
	\bigcap_{n\ge 0}\mathcal{F}^Y_+\vee\mathcal{F}^X_{[n,\infty[} =
	\mathcal{F}^Y_+
	\qquad\mathbf{P}^\nu\mbox{-a.s.}
$$
For the significance of this identity, we refer to \cite{Van08c} and the 
references therein.
\end{rem}

\section{Proof of Proposition \ref{prop:ar}}

We begin by proving the following representation.

\begin{lem}
For fixed $x,x'\in E$, we have
$$
	\|P(x,\cdot)-P(x',\cdot)\|_{\rm TV} =
	\int \left|
	\frac{q_\eta(\sigma(x)^{-1}\{\sigma(x')z-b(x)+b(x')\})}{
	\mathrm{det}(\sigma(x')^{-1}\sigma(x))}
	- q_\eta(z)
	\right|dz.
$$
\end{lem}

\begin{proof}
Note that $\sigma(x)$ is invertible for every $x$ as it is presumed to be
lower bounded.  Therefore $P(x,\cdot)$ has density
$p(x,z)=q_\eta(\sigma(x)^{-1}\{z-b(x)\})/\mathrm{det}(\sigma(x))$ with
respect to the Lebesgue measure on $\mathbb{R}^m$ for every $x$.
This implies that
$$
	\|P(x,\cdot)-P(x',\cdot)\|_{\rm TV} =
	\int \left|
	\frac{q_\eta(\sigma(x)^{-1}\{z-b(x)\})}{
	\mathrm{det}(\sigma(x))} -
	\frac{q_\eta(\sigma(x')^{-1}\{z-b(x')\})}{
	\mathrm{det}(\sigma(x'))}
	\right|dz.
$$
The result follows through a change of variables.
\end{proof}

We now prove that assumption \ref{aspt2} holds in this setting.

\begin{lem}
$\|P(x_n,\cdot)-P(y_n,\cdot)\|_{\rm TV}\to 0$ whenever $d(x_n,y_n)\to 0$.
\end{lem}

\begin{proof}
Fix any sequence $x_n,y_n$ such that $d(x_n,y_n)\to 0$.
By the previous lemma, it evidently suffices to show that the 
following function converges to $q_\eta(z)$ in $L^1(dz)$ as $n\to\infty$:
$$
	g_n(z) = 
	q_\eta(\sigma(x_n)^{-1}\{\sigma(y_n)z-b(x_n)+b(y_n)\})\,\mathrm{det}(
	\sigma(x_n)^{-1}\sigma(y_n)).
$$ 
Suppose first that $q_\eta$ is continuous.  Note that
$$
	\|\sigma(x_n)^{-1}\sigma(y_n)-I\| =
	\|\sigma(x_n)^{-1}\{\sigma(y_n)-\sigma(x_n)\}\|
	\le \alpha^{-1}\|\sigma(y_n)-\sigma(x_n)\|\xrightarrow{n\to\infty}0
$$
as $\sigma$ is uniformly continuous and lower bounded by $\alpha>0$, while
$\|\sigma(x_n)^{-1}\{b(x_n)-b(y_n)\}\|\le\alpha^{-1}\|b(x_n)-b(y_n)\|\to 
0$ as $n\to\infty$ as $b$ is uniformly continuous.  Therefore, if $q_\eta$ 
is continuous, then $g_n(z)$ converges to $q_\eta(z)$ pointwise.
By Scheff{\'e}'s lemma $g_n\to q_\eta$ in $L^1(dz)$.

Now suppose that $q_\eta$ is not continuous.  Then there is for every 
$\varepsilon>0$ a nonnegative continuous function with compact support 
$q_\eta^\varepsilon$ such that 
$\|q_\eta-q_\eta^\varepsilon\|_{L^1(dz)}<\varepsilon$ 
\cite[proposition 7.9]{Fol99}.  Using the triangle inequality, we easily 
estimate
$$
	\|P(x_n,\cdot)-P(y_n,\cdot)\|_{\rm TV} \le
	2\varepsilon +
	\int \left|
	\frac{q_\eta^\varepsilon(\sigma(x_n)^{-1}\{\sigma(y_n)z-b(x_n)+b(y_n)\})}{
	\mathrm{det}(\sigma(y_n)^{-1}\sigma(x_n))}
	- q_\eta^\varepsilon(z)
	\right|dz.
$$
But we have already established that the second term on the right 
converges to zero as $n\to\infty$, and $\varepsilon>0$ is arbitrary.  This 
completes the proof.
\end{proof}

Filter stability now follows from theorem \ref{thm:tv} whenever 
$\mathbf{P}^\mu|_{\sigma\{(Y_k)_{k\ge 0}\}}\ll
\mathbf{P}^\nu|_{\sigma\{(Y_k)_{k\ge 0}\}}$.  It remains to prove that 
when $q_\eta,q_\xi>0$ the absolute continuity requirement is in fact 
superfluous:
$$
	\|\pi_{n}^\mu-\pi_{n}^\nu\|_{\rm TV}
	\xrightarrow{n\to\infty}0\quad\mathbf{P}^\mu\mbox{-a.s.}
	\qquad\mbox{for any }
	\mu,\nu.
$$
Indeed, if this is the case, then by the triangle inequality
$$
	\|\pi_{n}^\mu-\pi_{n}^\nu\|_{\rm TV} \le
	\|\pi_n^\gamma-\pi_n^\mu\|_{\rm TV} +
	\|\pi_n^\gamma-\pi_n^\nu\|_{\rm TV}
	\xrightarrow{n\to\infty}0\quad
	\mathbf{P}^\gamma\mbox{-a.s.}
	\qquad\mbox{for any }\mu,\nu,\gamma,
$$
which completes the proof of proposition \ref{prop:ar}.

To establish the claim, note that when $q_\eta>0$ the transition kernel 
$P(x,\cdot)$ has a strictly positive density with respect to the Lebesgue 
measure for every $x\in E$.  In particular, $P(x,\cdot)\sim P(z,\cdot)$ 
for every $x,z\in E$.  As $q_\xi>0$ the filtering recursion is well 
defined under any initial measure, and it is immediately evident from the 
filtering recursion that $\pi_{1-}^\mu$ has a strictly positive density 
with respect to the Lebesgue measure for every $\mu$.  In particular, this 
implies that $\pi_{1-}^\mu\sim\pi_{1-}^\nu$ for every $\mu,\nu$.  But it 
is not difficult to establish that (see, e.g., the proof of \cite[lemma 
5.12]{Van08c})
$$
	\mathbf{E}^\mu\Big(\limsup_{n\to\infty}
	\|\pi_{n}^\mu-\pi_{n}^\nu\|_{\rm TV}\Big|Y_0=y\Big) = 
	\mathbf{E}^{\tilde\mu(y)}\Big(\limsup_{n\to\infty}
        \|\pi_{n}^{\tilde\mu(y)}-\pi_{n}^{\tilde\nu(y)}\|_{\rm TV}\Big),
$$
where $\tilde\mu(y)$ and $\tilde\nu(y)$ denote the regular conditional 
probabilities $\mathbf{P}^\mu(X_1\in\cdot|Y_0=y)$ and
$\mathbf{P}^\nu(X_1\in\cdot|Y_0=y)$ (i.e., $\tilde\mu(Y_0)=\pi^\mu_{1-}$ 
and $\tilde\nu(Y_0)=\pi^\nu_{1-}$).  As $\tilde\mu(y)\sim\tilde\nu(y)$
for (almost) every $y\in\mathbb{R}^n$, we have
$\mathbf{P}^{\tilde\mu(y)}|_{\sigma\{(Y_k)_{k\ge 0}\}}\sim
\mathbf{P}^{\tilde\nu(y)}|_{\sigma\{(Y_k)_{k\ge 0}\}}$ and the claim 
follows from theorem \ref{thm:tv}.

\begin{rem}
An almost identical argument shows that when $q_\eta,q_\xi>0$,
absolute continuity of the observations
$\mathbf{P}^\mu|_{\sigma\{(Y_k)_{k\ge 0}\}}\ll
\mathbf{P}^\nu|_{\sigma\{(Y_k)_{k\ge 0}\}}$ holds for any pair of initial 
measures $\mu,\nu$.  Together with theorem \ref{thm:tv} this gives the 
desired claim.  Note in particular that in this setting $\pi_n^\mu$ is 
$\mathbf{P}^\nu$-a.s.\ uniquely defined for any $\mu,\nu$, so that the 
statement of the proposition \ref{prop:ar} is in fact well posed 
(i.e., we do not need to be careful to choose a specific version of the 
filter).
\end{rem}

\bibliographystyle{amsplain}
\bibliography{ref}

\end{document}